\documentclass[11pt]{article}
\usepackage{epsfig,graphicx,amsfonts}

\usepackage{epsfig,theorem,latexsym,amsmath,subfigure}
\usepackage{epstopdf}

\def\dst{\displaystyle}

\def\Real{\mathbb R}

\def\dst{\displaystyle}

\def\bea{\begin{eqnarray}}
\def\eea{\end{eqnarray}}
\def\beann{\begin{eqnarray*}}
\def\eeann{\end{eqnarray*}}
\def\beeq#1{\begin{equation}{#1}\end{equation}}

\def\be{\begin{equation}}
\def\ee{\end{equation}}
\def\ba{\begin{array}}
\def\ea{\end{array}}
\def\bea{\begin{eqnarray}}
\def\eea{\end{eqnarray}}
\def\beann{\begin{eqnarray*}}
\def\eeann{\end{eqnarray*}}

\newtheorem{assumption}{Assumption}

\newtheorem{theorem}{Theorem}

\newtheorem{lemma}{Lemma}
\newtheorem{proposition}{Proposition}
\newtheorem{corollary}{Corollary}

\title{\LARGE \bf Robust output synchronization of a network of heterogeneous nonlinear agents via nonlinear regulation theory}
\author{A. Isidori, L. Marconi and G. Casadei
 \thanks{ A. Isidori is with DIS-``Sapienza" - Universit\`a di Roma, Roma, Italy {\tt \small albisidori@dis.uniroma1.it}}
 \thanks{ Lorenzo Marconi and G. Casadei are with C.A.SY. -- DEI,
 University of Bologna, Bologna, Italy  {\tt\small lorenzo.marconi@unibo.it}, {\tt\small g.casadei@unibo.it}.
 }
 }
   
\begin{document}

 \maketitle
\thispagestyle{empty}
\pagestyle{empty}

\begin{abstract}
In this paper, we consider the output synchronization problem for a network of heterogeneous diffusively-coupled nonlinear agents. Specifically, we show how the (non-identical) agents can be controlled in such a way that their outputs asymptotically track the output of a prescribed nonlinear exosystem. The problem is solved in two steps. In the first step, the problem of achieving consensus among (identical) nonlinear reference generators is addressed. In this respect, it is shown how the techniques recently developed to solve the  consensus problem among linear agents can be extended to agents modeled by nonlinear $d$-dimensional differential equations, under the assumption that the communication graph is connected. In the second step, the theory of nonlinear output regulation is applied in a decentralized control mode, to force the output of each agent of the network to robustly track the (synchronized) output of a local reference model.
 \end{abstract}

\section{Introduction}
The problem of achieving consensus (among states or outputs) in a (homogeneous or heterogenous) network of systems has attracted a major attention in the last decade. An exhaustive coverage of the literature, which is beyond the scope of the present paper, can be found, e.g. in the recent dissertation \cite{Wieland-diss} and in all references cited therein. We limit ourselves to mention that the case of a network of linear systems connected through a time-invariant graph has been fully addressed in the papers \cite{Scardovi-Sepulchre}, \cite{Seo-Shim-Back}, \cite{Wiel-Sep-All}, \cite{Kim-Shim-Seo}, while the analysis of a network of linear systems connected through a time-varying graph reposes on a fundamental convergence result established in \cite{Moreau-a}, \cite{Moreau-b}. Major results concerning the consensus problem in a network of nonlinear systems can be found in \cite{Hale}, \cite{Qu-Chu}, \cite{MMaggiore}, \cite{Arcak}, \cite{Stan-Sep}. The purpose of this paper is to present further contributions to the problem of output synchronization in a network of heterogeneous diffusively-coupled nonlinear agents.

As shown in \cite{Wiel-Sep-All} for linear systems and in \cite{Wieland-diss} for nonlinear systems, if the outputs of the agents of a heterogenous network achieve consensus on a nontrivial trajectory, the trajectory in question is necessarily the output of some autonomous (linear or nonlinear, depending on the case) system. This is the equivalent, in the context of the consensus problem, of the celebrated {\em internal model principle} of control theory \cite{FW76}. Motivated by this, we consider in what follows the problem of controlling a set of networked (non-identical) {\em nonlinear agents}  in such a way that their outputs asymptotically track the output of a prescribed {\em nonlinear} exosystem. The problem is solved in two steps.  In the first step, a network of $N$ identical copies of the given nonlinear exosystem is considered, the $k$-th of which is to be seen as ``local reference generator" for the $k$-th agent, and it is shown how certain  ``coupling gains" can be chosen in such a way that these local generators synchronize on a common consensus trajectory. To this end,  we extend existing techniques (see  \cite{Scardovi-Sepulchre}, \cite{Seo-Shim-Back}) recently proposed for the synchronization of a homogeneous network of linear systems  exchanging information through a connected (time-invariant) communication graph. The arguments used in this part are inspired by the literature on high-gain stabilization of nonlinear systems by output feedback, and specifically by the design of high-gain observers (see e.g. \cite{GK}). In the second step, the problem of controlling the  individual agent in such a way that its output tracks a reference output generated by the ``local exosystem" is addressed as a ``classical"  problem of nonlinear output regulation. In this respect, it is shown how the theory of nonlinear output regulation proposed in \cite{SICON} (see also \cite{BI2004}, \cite{CASYBook}),  can be successfully adopted to design robust internal model-based local regulators if the dynamics of the local agents fulfill a (weak) minimum-phase assumption. The result presented in the paper can be seen as a kind of  ``separation principle", in which the tools used to design local regulators  having the internal model property with respect to a local nonlinear exosystem in steady state, and the tools for synchronizing a set of networked nonlinear homogenous exosystems to reach  a common steady state, can be combined to achieve consensus of the outputs of heterogenous networked nonlinear systems.

The paper is organized as follow. In the next Section the problem is precisely formulated and the structure of the controllers is specified. Section III presents the solution for the first of the steps mentioned above, namely consensus in a network of  diffusively-coupled identical nonlinear systems. The problem of reaching a consensus  between heterogenous systems by means of nonlinear internal model-based regulator is addressed in Section IV, while Section V presents some simulation results, concerning the theory presented in Section III.

\section{Problem statement}

\subsection{Communication graphs.}

In what follows, the communication between individual systems (agents) is encoded by  a time-invariant {\em communication graph}. The latter is a triplet ${\mathcal G}=\{{\mathcal V}, {\mathcal E}, A\}$ in which:
\begin{itemize}
\item ${\mathcal V}$ is a set of $N$ {\em vertices} ${\mathcal V}=\{v_1,v_2, \ldots, v_N$\}, one for each of the $N$ agents in the set.
\item ${\mathcal E}\subset {\mathcal V}\times {\mathcal V}$ is a set of {\em egdes} that models the interconnection between nodes, according to the following convention: $(v_k,v_j)$ belongs to ${\mathcal E}$ if there is a flow of  information from node $j$ to node $k$.
\item the flow of information from node $j$ to node $k$  is {\em weighted} by the $(k,j)$-th entry $a_{kj}$ of the {\em adiacency matrix} $A\in \Real^{N\times N}$.
    \end{itemize}

    It is assumed that there are no self-loops, i.e. that  $(v_k,v_k)\notin {\mathcal E}$.   The set of {\em neighbors} of node $v_k$  is the set ${\mathcal N}_k=\{v_j \in {\mathcal V}: a_{k,j}\ne 0\}$. A {\em path} from node $v_j$ to node $v_k$ is a sequence of  $r$ distinct nodes $\{v_{\ell_1},  \ldots, v_{\ell_r}\}$ with $v_{\ell_1}=v_j$ and $v_{\ell_r}=v_k$ such that $(v_{i+1},v_i)\in{\mathcal E}$. A graph ${\mathcal G}$ is said to be {\em connected} if there is a node $v$ such that, for any other node $v_k \in {\mathcal V}\setminus\{v\}$, there is a path from $v$ to $v_k$.

 In what follows, we will consider cases in which the information available for control purpose at the $k$-th agent at time $t$ has the form
 \beeq{\label{nunu}
 \nu_k = \sum_{j=1}^N a_{kj}\, (\zeta_j(t)-\zeta_k(t))\qquad k=1,\ldots, N
 }
 in which $\zeta_i$, for $i=1,\ldots,N$, is a measurement taken at agent $i$. Letting $L$ denote the so-called matrix {\em Laplacian matrix} of the graph, defined by
 \[\ba{rclr}
 \ell_{kj}(t) &=& - a_{kj}(t) &\mbox{for $k\ne j$}\\[2mm]
 \ell_{kj}(t) &=& \sum_{i=1}^N a_{ki}(t) &\mbox{for $k= j$}\,,\ea\]
 the expression (\ref{nunu}) can be re-written as
\beeq{\label{nuell}
 \nu_k = -\sum_{j=1}^N\, \ell_{kj} \zeta_j(t)\qquad k=1,\ldots, N\,.
 }
By definition, the diagonal entries of $L$ are non-negative, the off-diagonal elements are non-positive and, for each row, the sum of all elements on this row is zero. As a consequence,  the all-ones $N$-vector $1_N=
 {\rm col}(1,1,\ldots, 1)\,$ is an eigenvector of $L$, associated with the eigenvalue
$\lambda=0$ . Let the other (possibly nontrivial) $N-1$ eigenvalues of $L$ be denoted as $\lambda_2(L), \ldots, \lambda_N(L)$.

 \begin{theorem} \label{TMconnect} A time-invariant graph is connected if and only if  its Laplacian matrix $L$ has only one trivial eigenvalue $\lambda_1=0$ and all other eigenvalues $\lambda_2(L), \ldots, \lambda_N(L)$ have positive real parts. \end{theorem}

\subsection{Problem formulation}
 We consider in what follows the problem of inducing consensus between the outputs of $N$ non-identical nonlinear systems,  which exchange information through a communication graph ${\cal G}$. The control system is decentralized, i.e. there is no leader sending information to each individual system, but rather each system exchanges information only with a set of neighboring systems, the information in question concerning only the {\em relative} values of the respective controlled outputs. { The $N$ nonlinear agents  are described by
 \beeq{\label{agents}
  \ba{rcl}
  \dot x_k &=& f_k(x_k) + g_k(x_k)u_k\\
  y_k &=& h_k(x_k)
  \ea
   \quad  x_k \in \Real^{n_k}\,, u_k, y_k \in \Real
 }
 $k = 1,\ldots, N$, where $u_k$ and $y_k$ are the local control input and output, with the inputs $u_k$ that must be designed in such a way that the
 outputs  $y_k$ of the $N$ systems asymptotically reach consensus on a nontrivial common trajectory $y^\ast(t)$.
 Each agent  is controlled by a local output-feedback controller of the form
 \beeq{\label{localcontroller}
 \ba{rcll}
  \dot \xi_k &=& \varphi_k(\xi_k, y_k,\nu_k) & \xi_k \in \Real^{\bar n_k}, \;\nu_k \in \Real^p\\
  u_k &=& \gamma_k(\xi_k,y_k,\nu_k) & u_k\in \Real\\
  \zeta_k &=& \rho_k(\xi_k,y_k) & \zeta_k \in \Real^p
  \ea
 }
 in which  $\zeta_k$ and $\nu_k$ are outputs and inputs that characterize the exchange of relative information between individual (controlled) agents, which takes the form (\ref{nunu}).

 In general terms, the control problem can be formulated as follows. Let $X_k\in \Real^{n_k}$, $k=1,\ldots, N$, be fixed compact set of admissible conditions for (\ref{agents}).
The problem is to find $N$ local controllers of the form (\ref{localcontroller}), exchanging information as in (\ref{nunu}), and compact sets $\Xi_k \in \Real^{\bar n_k}$, $k=1,\ldots, N$, of admissible initial conditions for all such controllers,  so that the positive orbit of the set of all admissible initial conditions is bounded and output consensus is reached, i.e. for each admissible initial condition $(x_k(0),\xi_k(0))\in X_k\times \Xi_k$, $k=1,\ldots, N$, there is a function $y^\ast : \Real \to \Real^d$  such that
\[
\lim_{t\to \infty}|y_k(t) - y^\ast(t)|=0 \qquad \mbox{for all $k=1,\ldots, N$}\,,
\]
uniformly in the initial conditions.

With the results of \cite{Wieland-diss} in mind, we expect that the consensus trajectory $y^\ast(t)$  can be thought of as generated by a {\em nonlinear} autonomous system, which could be modeled as an ordinary differential equation of order $d$
\beeq{\label{ystar}
 y^{\ast (d)} = \phi(y^{\ast
}, y^{\ast (1)}, \ldots, y^{\ast (d-1)})
 }
or in the equivalent state-space form of a $d$-dimensional system with output   \beeq{ \label{exo1}
   \ba{rcl}
  \dot w &=& s(w)\,, \qquad w \in \Real^d\\
  y^\ast &=& \vartheta(w)
  \ea
 }in which
 \beeq{\label{exo2}
 s(w) = Sw + B\phi(w)\,, \qquad \vartheta(w)=Cw
 }
and $(S,B,C)$ is a triplet of matrices in {\em prime} form. {Since we are seeking nontrivial consensus trajectories, in what follows we will consider the case in which (\ref{exo1}) possesses a nontrivial compact invariant set $W$. Moreover, we will assume that the function $\phi(\cdot)$ is globally Lipschitz.  In presence of systems of the form (\ref{exo2}) in which the $\phi(\cdot)$ is only locally Lipschitz, this assumption can be always enforced by properly modifying the function outside the compact set $W$ by using appropriate extension theorems.
}

  \subsection{Structure of local controllers and communication protocol}

 Bearing in mind the possibility of modeling all solutions of (\ref{ystar}) as outputs of the autonomous system (\ref{exo1})--(\ref{exo2}), in what follows, we consider for the local controllers (\ref{localcontroller}) a structure of the form
 \beeq{\label{bigstructure}
 \ba{rcl} \dot w_k &=& s(w_k) + \displaystyle K\sum_{j=1}^N a_{kj}(\vartheta(w_j) - \vartheta(w_k))\\
 \dot \eta_k &=& \varphi_k(\eta_k, e_k)\\
u_k &=& \gamma_k(\eta_k, e_k)\ea}
in which
\beeq{\label{errk}
 e_k = y_k - \vartheta(w_k)\,.
}

It is readily seen that this structure consists of a set of $N$ {\em local reference generators}
\beeq{\label{localexo}\ba{lcl}
\dot w_k &=& s(w_k) + K\nu_k\\
y^{\;{\rm ref}}_k &=& \vartheta(w_k)\,,\ea}
coupled via
\beeq{\label{exocoupl}
\nu_k=
\sum_{j=1}^N a_{kj}(y^{\;{\rm ref}}_j - y^{\;{\rm ref}}_k)\,,
}
each one of which provides a reference $y^{\;{\rm ref}}_k$ to be tracked by a {\em local regulator}
\[\ba{rcl}
\dot \eta_k &=& \varphi_k(\eta_k, e_k)\\
u_k &=& \gamma_k(\eta_k, e_k)\ea\]
driven by the local tracking error \[
e_k = y_k -y^{\;{\rm ref}}_k\,.\]

This control structure enables us to solve the problem in two stages. In the first stage, the design parameter $K$ is chosen in such a way as to induce consensus among the $N$ local generators (\ref{localexo}). In the second stage, the local regulators are designed in such a way that each of the outputs $y_k$ tracks its own reference $y^{\;{\rm ref}}_k$. It goes without saying that in the second step will ought to be able to use -- off the shelf  -- a large amount of existing results about the design of output regulators for nonlinear systems in the presence of exogenous signals generated by a nonlinear exosystem.

In this framework, we address first the problem of achieving consensus among the $N$ local generators (\ref{localexo}). Similar problems have received a large attention in recent literature (an exhausting covering of all such literature is beyond the scope of this paper, excellent surveys can be found in \cite{Wieland-diss}, \cite{Scardovi-Sepulchre}, \cite{Seo-Shim-Back}) but, to the best of our knowledge, a solution to the problem in the general terms considered here, i.e. for a network of $N$ {\em nonlinear} systems of dimension $d>1$, with information exchange in terms of relative values of $1$-dimensional outputs (other than relative values of their $d$-dimensional states), has not been proposed yet. We will address this problem in the following section.

\section{Achieving consensus in a homogeneous network of nonlinear systems}

As anticipated, we consider in what follows the problem of achieving {\em state} consensus in a network of $N$ identical nonlinear systems of the form (\ref{localexo}) coupled as in (\ref{exocoupl}), in which $s(w)$ and $\vartheta(w)$ are the map and the function defined in (\ref{exo2}).

The problem will be solved under the following assumptions.

  \begin{assumption}\label{ass1}
   The graph ${\cal G}$ is connected.
   \end{assumption}

 \begin{assumption}\label{ass2}
There exists a compact set  $W\subset \Real^d$ invariant  for (\ref{exo1}) such that the system
\[
  \dot w = S w + B \phi(w) + v
\]
 is input-to-state stable with respect to $v$ relative to $W$, namely there exist
a class-${\cal K}{\cal L}$ function $\beta(\cdot,\cdot)$ and a class-$\cal K$ function $\gamma(\cdot)$ such that
 \[
 \|w(t,\bar w)\|_{W} \leq \max\{ \beta(\|\bar w\|_W,t),\; \gamma(\sup_{\tau \in [0,t)}\|v(\tau)\|)\}\,.
 \]
 \end{assumption}

  \medskip
 To the purpose of inducing consensus in the network (\ref{localexo})--(\ref{exocoupl}), we choose the vector $K$ in (\ref{localexo}) as
 \[K=D_gK_0\,,\] where $D_g= \mbox{diag}\left( \ba{cccc} g & g^2 & \ldots & g^d \ea \right )$, with $g$  a design parameter and $K_0$  a vector to be designed.

 Due to Assumption 1, it is known that the Laplacian matrix $L$ has only one trivial eigenvalue and all remaining eigenvalues have positive real part. Hence there exists a $\mu>0$ such that
  \beeq{\label{Reeig}
   \mbox{Re}[\lambda_i(L)] \geq \mu \qquad i=2, \ldots, N\,.
  }

  With this in mind, let $T \in \Real^{N \times N}$ be defined as
  \[
   T= \left( \ba{cc} 1 & 0_{1 \times N-1}\\
   1_{N-1} & I_{N-1}
   \ea \right )
  \]
  and note that
  \[
   \tilde L = T^{-1} L T = \left( \ba{cc} 0 & L_{21}\\ 0_{N-1 \times 1} & L_{22} \ea \right )
  \]
  in which the eigenvalues of $L_{22}$ coincide with $\lambda_2(L), \ldots, \lambda_N(L)$. Then, the following result holds.

  \medskip
 \begin{lemma}Let $P$ be the unique positive definite symmetric solution of the algebraic Riccati equation
  \[
   S P + PS^T - 2 \mu P C^T C P + a I=0
  \]
  with $a>0$. Take $K_0$ as
  \beeq{\label{choiceK}
  K_0=P C^T\,.}
  Then, the matrix
  \[
  \left[ (I_{N-1} \otimes S) - (L_{22} \otimes K_0 C) \right ]
  \]
  is Hurwitz. $\triangleleft$
 \end{lemma}

  \medskip
 The proof of this Lemma can be found in \cite{Seo-Shim-Back} or in \cite{Wieland-diss}. Using this, we can now proceed with the proof of the main result of this Section.

  \medskip
 \begin{proposition}\label{prop1}
 Suppose Assumptions \ref{ass1} and \ref{ass2} hold. Consider the network of $N$ coupled systems
  \beeq{\label{exoscoupled}
\dot w_k = Sw_k + B\phi(w_k) + D_gK_0\sum_{j=1}^N a_{kj}(Cw_j - Cw_k)}
with $k=1,\ldots, N$. Let $K_0$ be chosen as in (\ref{choiceK}).
Then, there exists a number $g^\ast>0$  such that, for all $g \geq g^\ast$, the compact  invariant set
\beeq{\label{bigw}\ba{l}
{\bf W} = \{(w_1,w_2,\ldots,w_N)\in W\times W \times \cdots \times W : \\
\hspace*{2cm} w_1=w_2=\cdots= w_N\}\ea
}
is globally asymptotically stable. $\triangleleft$
\end{proposition}

 \medskip

\begin{corollary}
Let the hypotheses of the previous Proposition hold and let $K_0$ be chosen as in (\ref{choiceK}). There is a number $g^\ast>0$ such that, if $g\ge g^\ast$, the states of the $N$ systems (\ref{exoscoupled}) reach consensus, i.e. for every $w_k(0)\in \Real^d$, $k=1,\ldots, N$, there is a function $w^\ast : \Real \to \Real^d$  such that
\[
\lim_{t\to \infty}|w_k(t) - w^\ast(t)|=0 \qquad \mbox{for all $k=1,\ldots, N$}\,.\;\;\triangleleft
\]

 \end{corollary}
 \vspace*{2mm}
 {\bf Proof.}
 By the definition of Laplacian, the $k$-th controlled agent of the network (\ref{exoscoupled}) can be written as
 \[
 \dot w_k = Sw_k + B\phi(w_k) - D_g K_0 \sum_{j=1}^N  \ell_{kj} C w_j\,.
 \]
 Thus, setting ${\bf w}= {\rm col}(w_1,\ldots, w_N)$ entire set of $N$ controlled agents can be rewritten as
 \[
  \dot {\bf w} = \left[ (I_N \otimes S) - (L \otimes D_g K_0 C) \right ] {\bf w} + (I_N \otimes B) \Phi({\bf w})
 \]
 where
 \[
  \Phi({\bf w}) = \mbox{col}(\phi(w_1), \ldots, \phi(w_N))\,.
  \]

  Consider the change of variables ${\bf \tilde w} =( {T^{-1} \otimes I_d} ){\bf w}$, in which $T$ is the matrix introduced above. The system in the new coordinates read as
  \[
  \ba{rcl}
  \dot {\bf \tilde w} &=& 		
  (T^{-1} \otimes I_d) \left[ (I_N \otimes S) - (L \otimes D_g K_0 C) \right ]
  (T \otimes I_d) {\bf \tilde w}\\[1mm]
  && \hspace*{1cm} +\;(T^{-1} \otimes I_d) (I_N \otimes B) \Phi((T \otimes I_d){\bf \tilde w})\\
  &=&
  \left[ (I_N \otimes S) - (\tilde L \otimes D_g K_0 C) \right ] {\bf \tilde w} \\
  &&
  \hspace*{1cm}+\;(T^{-1} \otimes B) \Phi((T \otimes I_d){\bf \tilde w})\,.
  \ea
  \]
  Observing that \[
  {\bf \tilde w} = {\rm col}(w_1, w_2-w_1, \ldots, w_N-w_1)\] set $z_k= w_k-w_1$, for $k=2,3, \ldots, N$, and \[
  {\bf z} ={\rm col}(z_2,z_3, \ldots, z_N)\,,\]
  yielding ${\bf \tilde w}= {\rm col}(w_1,{\bf z})$.
  Then, it is readily seen that the system above
exhibits a triangular structure of the form
  \[
  \ba{rcl}
   \dot w_1 &=& S w_1 + B \phi(w_1) - (L_{12} \otimes D_gK_0 C) {\bf z}\\[2mm]
   \dot {\bf z} &=& \left[ (I_{N-1} \otimes S) - (L_{22} \otimes D_g K_0 C) \right ] {\bf z}
   +\Delta \Phi(w_1, {\bf z})
   \ea
  \]
  where
  \[
   \Delta \Phi(w_1, {\bf z}) = (I_{N-1} \otimes B) \left( \ba{ccc} \phi(w_1+z_2) - \phi(w_1)\\
    \phi(w_1+z_3) - \phi(w_1)\\
    \cdots\\
    \phi(w_1+z_N) - \phi(w_1)
   \ea \right )\,.
  \]
 Note that $\Delta \Phi(w_1, {\bf z})$ is globally Lipschitz in ${\bf z}$ uniformly in $w_1$ and $\Delta \Phi(w_1, 0) \equiv 0$ for all $w_1 \in \Real^n$.
 Consider now the rescaled state variable
 \[
 \zeta = (I_{N-1} \otimes D_g^{-1}) {\bf z}\,.
 \]
 By using the definition of $S$, $C$ and $B$, it follows that

  \[
  \ba{rcl}
   \dot w_1 &=& S w_1 + B \phi(w_1) - (L_{12} \otimes D_g K_0 C) (I_{N-1} \otimes D_g) \zeta\\[2mm]
   \dot \zeta &=&g  \left[ (I_{N-1} \otimes S) - (L_{22} \otimes K_0 C) \right ] \zeta\\[2mm]
   && \quad
   \hspace*{1cm}+\;\dst {1 \over g^{d}}\Delta \Phi(w_1, (I_{N-1} \otimes D_g) \zeta)\,.
   \ea
  \]
 It is known from the Lemma above that the proposed choice of $K_0$ guarantees that the matrix
 $ \left[ (I_{N-1} \otimes S) - (L_{22} \otimes K_0 C) \right ]$ is Hurwitz. As a consequence standard high-gain arguments lead to the conclusion that
 if $g$ is chosen sufficiently large, the equilibrium $\zeta=0$ of the lower subsystem is asymptotically stable, uniformly in $w_1$, and actually with a quadratic Lyapunov function that is independent of $w_1$. The ISS property in Assumption 2 thus guarantees that $w_1$ converges to the invariant set $W$. Since  $w_k = w_1 + z_k$, $k=2,\ldots,N$, and $z_k \to 0$ as $t\to \infty$, the result follows. $\triangleleft$

\section{Achieving consensus in a heterogenous network of nonlinear systems}

We proceed now with the second step of the design, i.e. we design local regulators for each agent. In what follows, we assume that the vector $K_0$ and the values of $g$ have been fixed and such that the conclusion of Proposition \ref{prop1} holds, i.e. such that the set (\ref{bigw}) is globally asymptotically stable for (\ref{exoscoupled}).

In what follows, we assume that each individual agent has a well defined relative degree $r$ between input $u_k$ and output $y_k$ and possess a globally defined normal form. To streamline the exposition, we consider the special case in which $r=1$. The case of higher relative degree only entails heavier notational complexity and no conceptual differences. Thus we assume that the individual agent is modeled by equations of the form
\beeq{\label{agentnf}
\ba{rcl}
\dot z_k &=& f_k(z_k,y_k)\\
\dot y_k &=& q_k(z_k,y_k) + b_k(z_k,y_k)u_k\,, \ea
}
where $z_k\in \Real^{n_k-1}$ and where $b_k(z_k, y_k)$, which is the high-frequency gain of the k-th agent, is bounded away from zero. In particular, we assume that there exists a $\bar b_k >0$ such that $b_k(z_k, y_k) > \bar b_k$ for all $(z_k,y_k) \in \Real^{n_k-1} \times \Real$ and for all $k=1,\ldots,N$.  As anticipated, with this system we associate a local tracking error of the form
\[
e_k = y_k-C w_k\,.
\]
The problem is to design a local regulator to the purpose of steering $e_k$ to zero.

In this respect, it should be borne in mind that $w_k$ is a ``portion" of the state of the coupled system (\ref{exoscoupled}) and hence the entire dynamics of the latter should be taken into account in the analysis. To this end,
invoking the arguments used in the proof of Proposition \ref{prop1},  observe that, for any $k$, there exists a linear change coordinates in (\ref{exoscoupled}),  by means of which  this system (the entire set of $N$ networked local reference generators) can be changed into a system modeled by equations of the form
\beeq{\label{superexok}
\ba{rcl}
\dot \zeta_k &=& \psi(\zeta_k,w_k)\\ \dot w_k &=& s(w_k) + \Upsilon_k \zeta_k
\ea
}
in which $\psi(0,w_k)=0$. By assumption, the lower sub-subsystem is input-to-state stable, with respect to the input $\zeta_k$, to the set $W$. { Moreover, as observed in the proof of Proposition \ref{prop1}, the equilibrium $\zeta_k=0$ of the upper subsystem is globally exponentially stable. As a consequence, the set
\[
 \{(\zeta_k,w_k)\in \Real^{(N-1)d}\times \Real^d: \zeta_k=0, w_k\in W\}
\]
is a globally asymptotically stable compact invariant set of (\ref{superexok}).}
%FINE secondo tratto rosso

In view of this, we can represent the aggregate of (\ref{agentnf}) and of (\ref{superexok}) as a standard ``exosystem-plant" interconnection
\beeq{\label{complete0}
\ba{rcl} \dot \zeta_k &=& \psi(\zeta_k,w_k)
\\ \dot w_k &=& s(w_k) + \Upsilon_k \zeta_k\\
\dot z_k &=& f_k(z_k,y_k)\\
\dot y_k &=& q_k(z_k,y_k) + b_k(z_k,y_k)u_k\\
e_k &=& y_k-Cw_k\,. \ea
}
As usual, we change variables replacing $y_k$ by $e_k$ and obtain
\beeq{\label{complete}
\ba{rcl} \dot \zeta_k &=& \psi(\zeta_k,w_k)
\\ \dot w_k &=& s(w_k) + \Upsilon_k \zeta_k\\
\dot z_k &=& f_k(z_k,Cw_k+e_k)\\
\dot e_k &=& q_k(z_k,Cw_k+e_k)- C[s(w_k)
+\Upsilon_k \zeta_k]  \\ &&\qquad+\; b_k(z_k,Cw_k+e_k)u_k \,.\ea
}
This system is ready for the design (under appropriate hypothesis) of a local regulator of the form
\beeq{ \label{RegGen}
\ba{rcl}\dot \eta_k &=& \varphi_k(\eta_k) +G_k v_k\\
u_k &=& \gamma_k(\eta_k) + v_k\\
v_k &=& \kappa_k(e_k)
\ea
}
according to the procedures suggested in \cite{BI2004} or in \cite{SICON}. The basic assumption needed to make the design possible is that the zero dynamics of (\ref{complete}) namely, those of
\beeq{\label{zerodynbig}
\ba{rcl} \dot \zeta_k &=& \psi(\zeta_k,w_k)
\\ \dot w_k &=& s(w_k) + \Upsilon_k \zeta_k\\
\dot z_k &=& f_k(z_k,Cw_k)\ea
}
possess a compact invariant set which is asymptotically stable with a domain of attraction that contains the prescribed set of initial conditions. { To make this assumption precise, let $W_k$ be the set of admissible conditions of $w_k$, let $S_k$ be the set of admissible initial conditions of $\zeta_k$ and $Z_k$  the set of admissible initial conditions of $z_k$. Then, the assumption is question can be stated as follows.

\medskip
\begin{assumption}\label{ass3}
There exist a (possibly set-valued) map $\pi_k: w_k\in W \mapsto \pi_k(w_k) \subset \Real^{n_k-1}$ such that the set
\[\ba{rcl}
{\mathcal A}_k &=& \{(\zeta_k,w_k,z_k):\zeta_k=0, w_k\in W, z_k\in \pi_k(w_k)\}\\ [1mm]
\ea
\]
is an asymptotically stable invariant set for (\ref{zerodynbig}) with a domain of attraction containing $S_k\times W_k\times Z_k$.
\end{assumption}

\medskip
We note that this assumption is the natural formulation, in the current framework of a networked system, of the (weak) minimum-phase assumption that one would assume in solving a problem of output regulation for the $k$-th agent if high-gain arguments were to be used for stabilization purposes.

We proceed now with the design the functions $(\varphi_k(\cdot), \gamma_k(\cdot), G_k)$  in (\ref{RegGen}), whose key properties are captured in the following definition, taken from  \cite{CASYBook}.\\[2mm]
{\em Definition} ({\em Asymptotic internal model property}).
The triplet $(\varphi_k(\cdot), \gamma_k(\cdot), G_k)$ is said to have the {\em asymptotic internal model property} if
 there exists a $C^1$ map $\tau_k: \Real^d \times \Real^{n_k-1} \to \Real^{m_k}$ such that the following holds: \footnote{We use the notation ${\rm gr}(\pi_k):=\{(w_k,z_k): w_k\in W, z_k \in \pi_k(w_k)\}$}

\smallskip\noindent (i)
for all $(w_k,z_k) \in {\rm gr}(\pi_k)$
 \[
  \ba{rcl}
   \dst {\partial \tau_k \over \partial w_k}s(w_k) + \dst {\partial \tau_k \over \partial z_k}f(z_k,Cw_k)&=& \varphi_k(\tau(w_k,z_k))\\[3mm]
   \dst {Cs(w_k) - q_k(z_k,Cw_k)  \over b_k(w_k,Cw_k)} &=& \gamma_k(\tau(w_k,z_k))
  \ea
 \]
\smallskip\noindent (ii)
  the set
\[\ba{l}
 {\mathcal S}_k = \{(\zeta_k,w_k,z_k,\eta_k): \\[1mm] \hspace*{1cm}
 \zeta_k=0, (w_k,z_k)\in {\rm gr}(\pi_k), \eta_k = \tau(w_k,z_k)\}\ea
 \]
is locally asymptotically stable for the system
\[
\ba{rcl}
\dot \zeta_k &=& \psi(\zeta_k,w_k)\\
\dot w_k &=& s(w_k) + \Upsilon_k \zeta_k\\
\dot z_k &=& f_k(z_k,Cw_k)\\
\dot \eta_k &=& \varphi_k(\eta_k) - G_k \Bigl[ \gamma_k(\eta_k) \\ [2mm]
&& \hspace*{1cm}+\, \dst {q_k(z_k,Cw_k) - C[s(w_k)+\Upsilon_k\zeta_k] \over b_k(w_k,Cw_k)}\Bigr]
 \ea
\]
with a domain of attraction containing $S_k\times W_k\times Z_k \times M_k$, where $M_k$ is the compact set of initial conditions of (\ref{RegGen}). $\triangleleft$

\medskip
If a triplet with the asymptotic internal model property can be designed then the problem of  steering the regulation error $e_k$ of the $k$-th agent to zero is solved as claimed by the following theorem proved in \cite{SICON}.}

 \medskip
\begin{theorem}
Let $S_k \subset \Real^{(N-1)d}$, $W_k\subset \Real^d$, $Z_k \subset \Real^{n_k-1}$, $E_k \subset \Real$ and $M_k\subset \Real^{m_k}$  be compact sets of initial conditions for the closed-loop  system (\ref{complete}), (\ref{RegGen}). Let the triplet $(\varphi_k(\cdot), \gamma_k(\cdot), G_k)$ be designed so that it has the asymptotic internal model property. Then there exists a continuous function $\kappa_k: \Real \to \Real$ such that the trajectories of the closed-loop system originating form $S_k\times W_k\times Z_k \times E_k \times M_k$ are bounded and $\lim_{t \to \infty} e_k(t)=0$ uniformly in the initial conditions.
\end{theorem}

\medskip
A triplet having the internal model property can be always designed as  detailed in the next result coming from a slight adaptation of the results presented in \cite{SICON}.

\medskip
\begin{proposition}
 Let $m_k \geq 2 (d+n_k-1) +2$. Then there exists a $\lambda_k^\ast<0$ and, for almost all possible choice of controllable pairs $(F_k,G_k)\in \Real^{m_k \times m_k} \times \Real^{m_k \times 1}$ such that $\mbox{Re}\sigma(F_k) \leq \lambda_k^\ast$, there exists a continuous $\gamma_k: \Real^{m_k} \to \Real$, such that the triplet $(\varphi_k(\cdot), \gamma_k(\cdot), G_k)$ with
 $\varphi_k (\eta_k)= F_k \eta_k + G_k \gamma_k(\eta_k)$ has the asymptotic internal model property.
\end{proposition}

\medskip
The previous results, although conceptually interesting, is not constructive in the design of the function $\gamma_k$. As shown in \cite{BI2004}, it turns out that a constructive design procedure can be given if an extra assumption is invoked. In particular, assume that there exists a $m_k>0$ and a locally Lipschitz function $\varrho_k: \Real^{m_k} \to \Real$ with the property that, for all $(w_k(0),z_k(0))\in {\rm gr}(\pi_k)$, the solution $w_k(t),z_k(t)$ of
\[\ba{rcl}
\dot w_k &=& s(w_k) \\
\dot z_k &=& f_k(z_k,Cw_k)
\ea\]
is such that the function \[\rho(t) =  \dst {Cs(w_k(t)) - q_k(z_k(t),Cw_k(t))  \over b_k(w_k(t),Cw_k(t))}\]
satisfies
\[
 \rho^{(m_k)}(t) = \varrho_k ( \rho(t),  \rho^{(1)}(t), \ldots,  \rho^{(m_k-1)}(t)) \quad \forall \, t \in \Real\,.
\]
If this assumption holds then the following result can be proved, by means of a slight adaptation of the results presented in \cite{CASYBook}.

\begin{proposition}
 Let $(A_k, B_k,C_k) \in \Real^{m_k \times m_k} \times \Real^{m_k \times 1} \times \Real^{1 \times m_k}$ be a triplet of matrices in prime form. Furthermore, let $\bar \varrho_k: \Real^{m_k} \to \Real$  a bounded locally Lipschitz function that agrees with $\varrho_k(\cdot)$ on $B_R=\{\xi \in \Real^{m_k}: |\xi|\le R\}$,
  let $D_\ell = \mbox{\rm diag}(\ell, \ell^2, \ldots, \ell^{m_k})$ with $\ell$ a positive design parameter, and let $(c_0, \ldots, c_{m_k-1})$ be such that the polynomial
 $\lambda^{m_k} + c_0 \lambda^{m_k-1} + \ldots + c_{m_k-1}$ is Hurwitz. Then there exist $R>0$ and $\ell^\ast>0$ such that for all $\ell \geq \ell^\ast$ the triplet $(\varphi_k(\cdot), \gamma_k(\cdot), G_k)$ defined as
 \[
  \varphi_k(\eta_k) = A_k \eta_k + B_k \bar \varrho_k(\eta_k)\,, \quad \gamma_k(\eta_k) = C_k \eta_k\,,
  \]
 $G_k = D_\ell \; \mbox{col}( c_0,\, \ldots, \,c_{m_k-1})$ has the asymptotic internal model property.
\end{proposition}

\section{Simulation results}
In this section we present simulation results about the theory presented in Section III by considering two different nonlinear oscillators as
system (\ref{exo1}), namely a Van der Pol and a Duffing oscillator. We consider the case of three agents ($N=3$).

In case of Van der Pol, system (\ref{exo1}), (\ref{exo2}) takes the form
\[
 \ba{rcl}
  \dot w_1 &=& w_2\\
  \dot w_2 &=& 2 (1 - w_1^2) w_2 - w_1
 \ea \qquad   y^\star = w_1\,.
\]
The topology of the simulated network is described by the incidence matrix
\[
 A= \left( \ba{ccc}
  0 & 0 & 1\\
   1 & 0 & 0\\
    0 & 1 & 0
 \ea
 \right )
\]
with the eigenvalues of the corresponding Laplacian that fulfill (\ref{Reeig})   with $\mu =1.4$.
Figures \ref{fig:VDP1}, \ref{fig:VDP2} present the  simulation results obtained with the three  systems of the form (\ref{exoscoupled}) where the coupling term has been fixed as in Lemma 1 with $g=10$ and $a=1$.  In particular,    Figure \ref{fig:VDP1} shows the phase-plane of the three oscillators initialized respectively at $w_1=(1, 1)$, $w_2=(3,3)$ and $w_3=(5,5)$, while Figure \ref{fig:VDP2} shows the time behaviors of their first state variable.

\begin{figure}[!htb]
\centering
\includegraphics[width=1\textwidth]{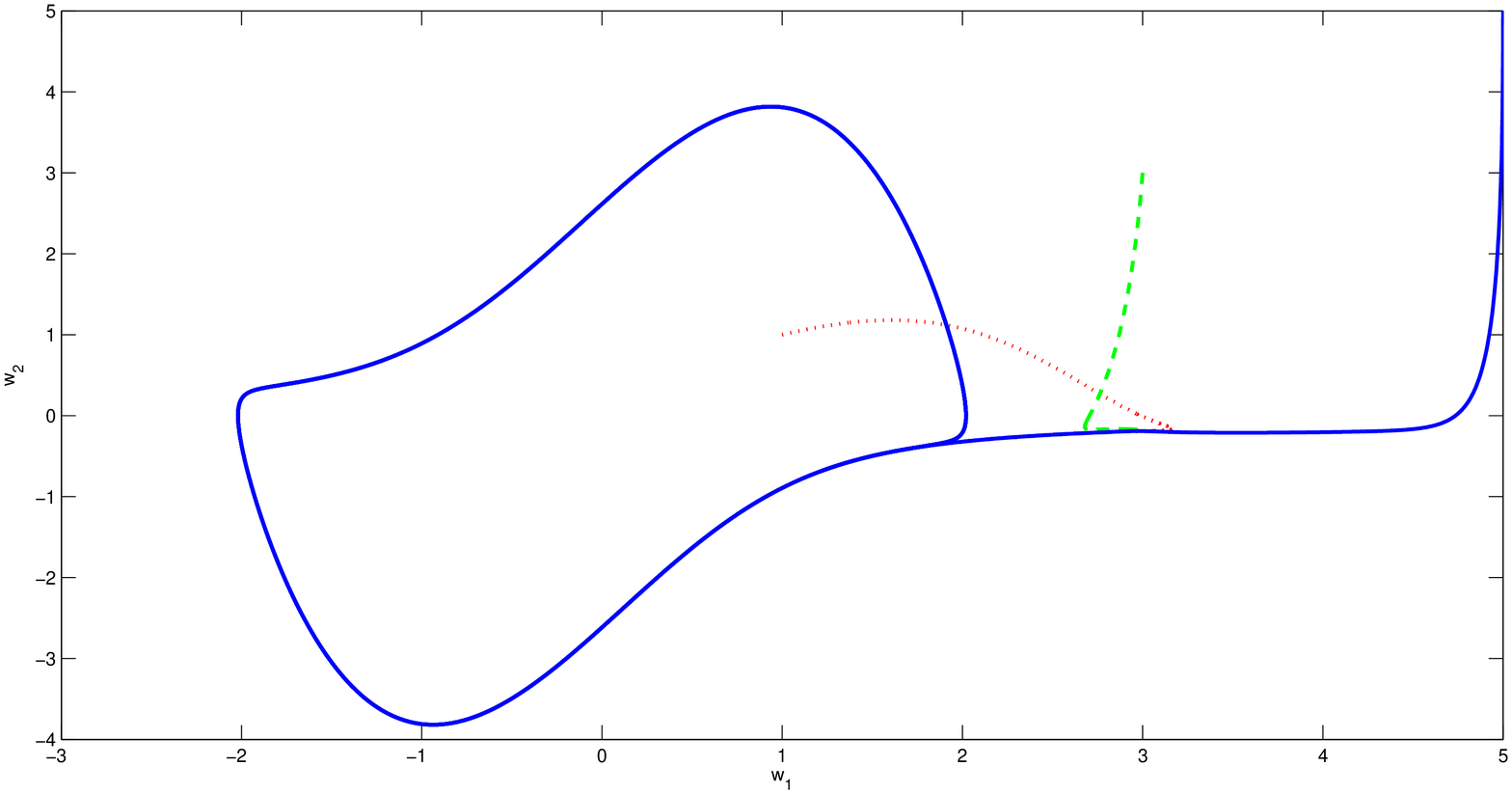}
\caption{\label{fig:VDP1} Phase-plane of the three networked Van der Pol Oscillators. }
\end{figure}

\begin{figure}[!htb]
\centering
\includegraphics[width=1\textwidth]{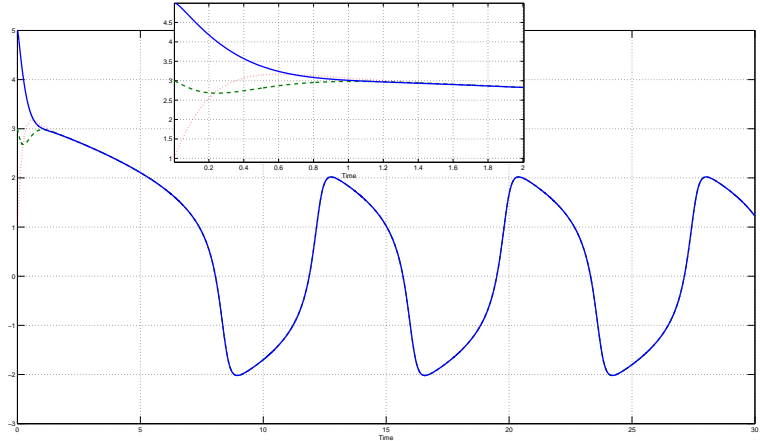}
\caption{\label{fig:VDP2} Time behavior of the first state variable of the three Van der Pol Oscillators. }
\end{figure}

As a second example we consider the case of a Duffing oscillator described by
\[
 \ba{rcl}
  \dot w_1 &=& w_2\\
  \dot w_2 &=&  -2 w_1  +  w_1^3
 \ea \qquad   y^\star = w_1\,.
\]
The topology of the network is described by the same incidence matrix considered in the previous example with the same
parameters that have been used in the coupling terms. Figure \ref{fig:DUFF1} shows the phase-plane of the three oscillators initialized respectively at $w_1=(1, 1)$, $w_2=(3,3)$ and $w_3=(5,5)$, while Figure \ref{fig:DUFF2} shows the time behaviors of their first state variable.

\begin{figure}[!htb]
\centering
\includegraphics[width=1\textwidth]{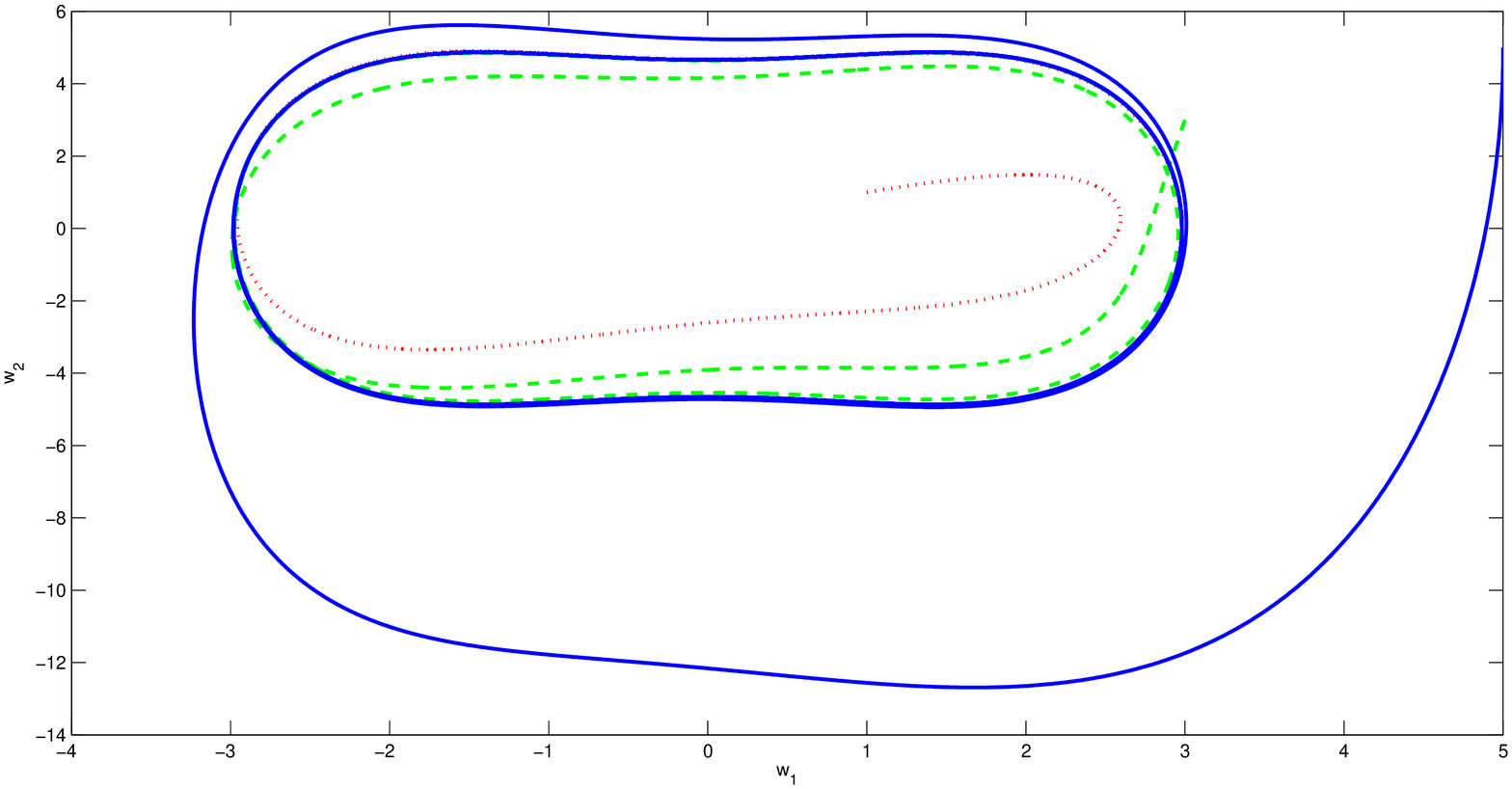}
\caption{\label{fig:DUFF1} Phase-plane of the three networked Duffing Oscillators. }
\end{figure}

\begin{figure}[!htb]
\centering
\includegraphics[width=1\textwidth]{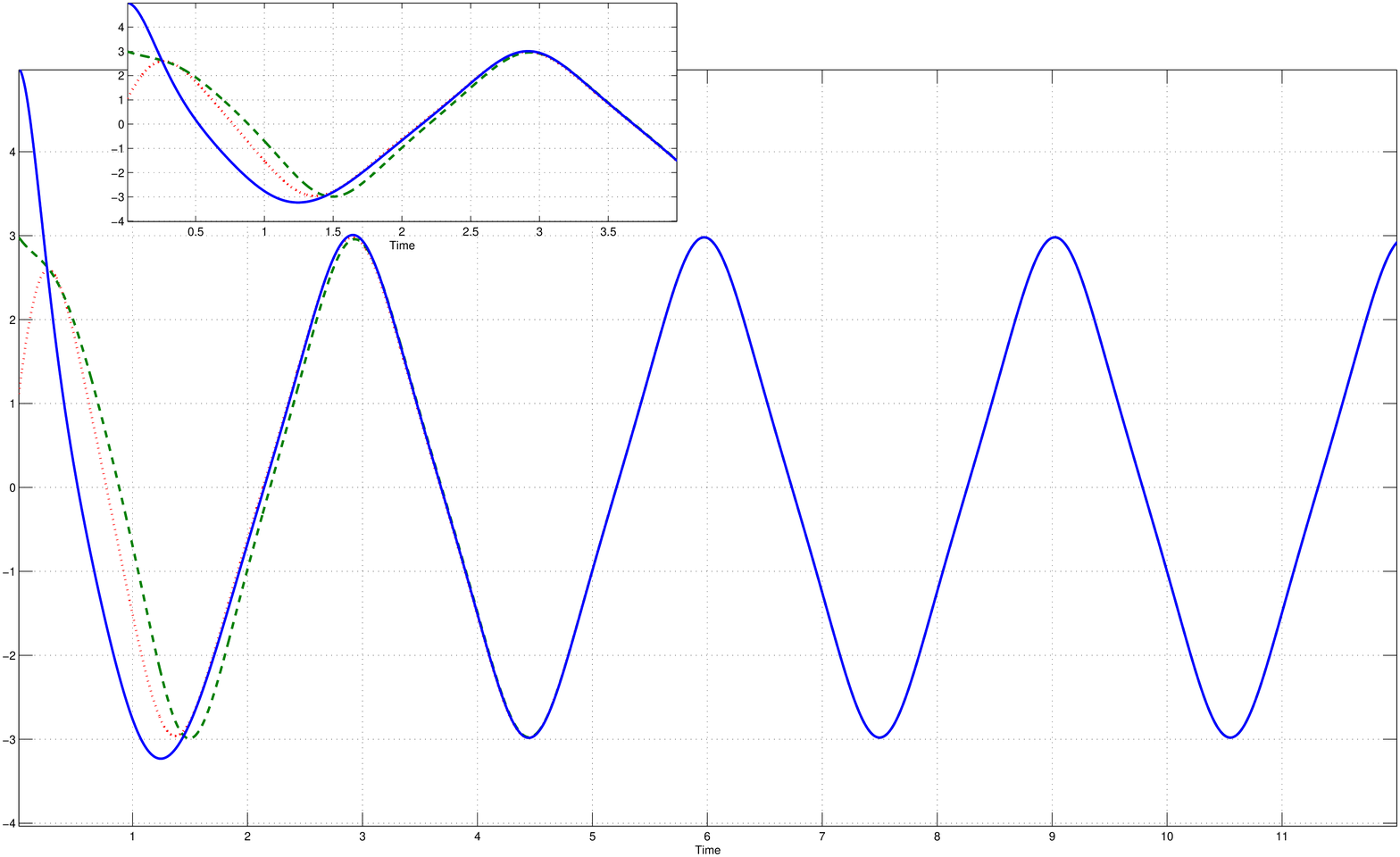}
\caption{\label{fig:DUFF2} Time behavior of the first state variable of the three Duffing Oscillators. }
\end{figure}

\section{Conclusions}

The problem of reaching a consensus among the outputs of a set of networked nonlinear systems was considered. The output reference signal is thought of as generated by the autonomous nonlinear exosystem of the form (\ref{ystar}). A first result of the paper is presented in Section III where it is shown how  $N$  diffusively-coupled exosystems of the form (\ref{localexo}) reach a consensus over a trajectory of (\ref{ystar}). This result is instrumental to the theory presented in Section IV where the nonlinear output regulation theory is adopted to design local regulators that make the outputs of $N$ heterogenous systems to track a common reference signals solution of  (\ref{ystar}).


\begin{thebibliography}{99}

\bibitem{Arcak} M. Arcak. Passivity as a design tool for group coordination. {\em IEEE Trans. Aut. Contr.}, 52(8), pp. 1380–1390, 2007.


\bibitem{BI2004} C. I. Byrnes, A. Isidori, ``Nonlinear internal models for output
regulation", {\em IEEE Trans. on Aut. Contr.}, 49(12), pp. 2244-2247, 2004

\bibitem{FW76}
B.A. Francis and W.M. Wonham,
\newblock ``The internal model principle of control theory'',
\newblock { \em Automatica}, 12, pp. 457--465, 1976.

\bibitem {GK}
J.P. Gauthier, I. Kupka, \newblock { \em Deterministic
Observation Theory and Applications}, \newblock Cambridge
University Press, Cambridge (2001).

\bibitem{Hale} J. K. Hale. ``Diffusive coupling, dissipation, and synchronization". {\em J. of Dynamics and
Diff. Eqns}, 9(1), pp. 1–52, 1997.



 \bibitem{Kim-Shim-Seo}   H. Kim, H.  Shim, and J.H. Seo. ``Output consensus of heterogeneous uncertain linear multi-agent
systems", {\em IEEE Trans. Aut. Contr.}, 56(1):  200-206, 2011.

\bibitem{MMaggiore} Z. Lin, B. Francis, M. Maggiore, ``State agreement for continuous-time coupled nonlinear systems", {\em SIAM J. Contr. Optimiz.} 46(1), pp. 288-307. 2007



\bibitem{CASYBook} L. Marconi and Isidori,  ``A unifying {approach to the design of} nonlinear output {regulators}". In {\em Advances in Control Theory and Applications}, C. Bonivento, A. Isidori, L. Marconi, C. Rossi Eds. , LNCIS, Springer Verlag Berlin, 2007
 \bibitem{SICON}  L. Marconi, L. Praly, A. Isidori, ``Output Stabilization via Nonlinear Luenberger Observers". {\em SIAM
J. on Contr. and Optimiz.}, 45(6),  pp. 2277-2298, 2007.


\bibitem{Moreau-a} L. Moreau ``Stability of continuous-time distributed consensus algorithms".
http://arxiv.org/abs/math/0409010v1,2004a. arXiv:math/0409010v1[math.OC].

\bibitem{Moreau-b} L. Moreau, ``Stability of multi-agent systems with time-dependent communication links", 50(2), pp. 169-182, 2005.

\bibitem{Qu-Chu} Z. Qu, J. Chunyu, and J.Wang, "Nonlinear cooperative control for consensus
of nonlinear and heterogeneous systems," in {\em Proc. IEEE Conf.
Decision Control},  pp. 2301–2308, 2007.


\bibitem{Scardovi-Sepulchre} L. Scardovi, R. Sepulchre, ``Synchronization in networks of identical linear systems", {\em Automatica}, 45(10), pp. 2557-2562, 2009.


\bibitem{Seo-Shim-Back} J. H. Seo,  H. Shima, J. Back, ``Consensus of high-order linear systems using dynamic output feedback compensator: low gain approach", {\em Automatica}, 45(11),   pp. 2659-2664, 2009.
    
 
\bibitem{Stan-Sep} G.-B. Stan and R. Sepulchre. ``Analysis of interconnected oscillators by dissipativity theory". {\em IEEE Trans. Aut. Contr.}, 52(2), pp. 256–270, 2007.


    
\bibitem{Wieland-diss} P. Wieland, {\em From static to Dynamic Couplings in Consensus and Synchronization among Identical and Non-Identical Systems}, PhD thesis, Universit\"at Stuttgart, 2010.

\bibitem{Wiel-Sep-All} P. Wieland, R.Sepulchre, and F. Allg\"{o}wer. ``An internal model principle is necessary and sufficient for linear output synchronization''. {\em Automatica}, 47, 1068-1074, 2011.

 \end{thebibliography}
\end{document}